\documentclass[sigconf, authorversion]{acmart} 

\AtBeginDocument{%
  }

\setcopyright{none}

\acmConference{none}{none}{none}
\acmISBN{}

\usepackage{xspace}
\usepackage{cleveref}
\usepackage{xcolor}
\usepackage{makecell}
\usepackage{subcaption}
\usepackage{adjustbox}
\usepackage{enumitem}
\usepackage{minted}

\setlist[enumerate]{label=\roman*)}
\setlist[itemize]{label=\mySqBullet}

\crefname{table}{Tab.}{Fig.}
\crefname{figure}{Fig.}{Fig.}
\crefname{section}{Sec.}{Sec.}
\crefname{equation}{Eq.}{Eq.}
\creflabelformat{equation}{#2\textup{#1}#3}
\usepackage{amsmath}
\usepackage{amsfonts}       

\definecolor{darkyellow}{rgb}{0.95, 0.73, 0.24}

\newcommand{\globe}{\textsc{GLOBe}\@\xspace}

\newcommand{\bb}[1]{\mathbb{#1}}

\newcommand{\dif}{\mathrm{d}}

\newcommand{\inlinetitle}[2] {\medskip\noindent\textbf{#1{#2}}}
\newcommand{\argmin}[1]{\underset{#1}{\arg \min} \;}
\newcommand{\argmax}[1]{\underset{#1}{\arg \max} \;}

\renewcommand{\epsilon}{\varepsilon}

\newcommand{\footurl}[1]{\footnote{\url{#1}}}

\usepackage{hyperref}
\definecolor{mydarkblue}{rgb}{0,0.08,0.45}
\definecolor{mylightpurple}{rgb}{0.26, 0.32, 0.82}
\definecolor{mylightblue}{rgb}{0.26, 0.38, 0.93}
\hypersetup{ %
  colorlinks=true,
  linkcolor=red,
  citecolor=mylightpurple,
  filecolor=mylightpurple,
  urlcolor=mylightpurple,
}

\usepackage[labelfont=bf]{caption}
\DeclareCaptionLabelSeparator{point}{.~}
\usepackage[normalem]{ulem}

\newcounter{marginNoteCounter}

\usepackage{xspace}
\newcommand{\eg}   			{e.g.\@\xspace}

\renewcommand{\ldots}{...}

\newcommand{\mySqBullet}		{\raisebox{0.25em}{{\scriptsize$_\blacksquare$}}}

\begin{document}
\title{\globe: A Modular Global Optimization library}


\author{Gaëtan Serré}
\email{name.surname@ens-paris-saclay.fr}
\affiliation{%
  \institution{École Normale Supérieure Paris-Saclay - Centre Borelli}
  \city{Gif-Sur-Yvette}
  \country{France}
}
\author{Argyris Kalogeratos}
\email{name.surname@ens-paris-saclay.fr}
\affiliation{%
  \institution{École Normale Supérieure Paris-Saclay - Centre Borelli}
  \city{Gif-Sur-Yvette}
  \country{France}
}
\author{Nicolas Vayatis}
\email{name.surname@ens-paris-saclay.fr}
\affiliation{%
  \institution{École Normale Supérieure Paris-Saclay - Centre Borelli}
  \city{Gif-Sur-Yvette}
  \country{France}
}

\renewcommand{\shortauthors}{Serré, Kalogeratos, and Vayatis}

\begin{abstract}
Open-source libraries are have a catalytic role in research pipelines, where new methods must be compared against up-to-date baselines. We present the GLobal Optimization Benchmark (\globe) modular Python library that unifies classical and recent continuous global optimization algorithms, including decision-based and mathematically founded particle-based methods, in a single framework. A central contribution of \globe is its modular architecture, which factors common algorithmic patterns into reusable family-level components and enables plugins to be implemented once and makes them available to all algorithms in the corresponding family. This modular design leverages recent advances in mathematical formalization of global optimization, where structural commonalities across algorithms have been identified and used to develop broadly applicable, formally grounded algorithmic features. The C++ backend relies on Eigen $5$ for efficient linear algebra. \globe presently includes $14$ optimizers, $19$ analytical benchmarks along with a random function generator, and an integrated toolkit for direct algorithm comparison.
\end{abstract}

\keywords{Global optimization, open-source, modular library, decision-based algorithms, particle-based algorithms, benchmarking.}


\maketitle

\section{Introduction}

Global optimization aims at finding the global optimizer (maximizer or minimizer) of a real-valued function $f$ over a search space $\Omega$, usually a subset of $\bb{R}^d$. 
The applications of global optimization are vast, ranging from machine learning and data science to operations research and engineering design \cite{Pinter1996, lee2017finding, Agrawal2021, Awasthi2024, Houssein2024, el2024taxonomy, Franceschi2024}. Although the study of global optimization algorithms dates decades back \cite{Kirkpatrick1983,Price1983}, it remains an active research area with many recent contributions, notably on SDE-driven methods \cite{Pinnau2017, Bungert2024, riedl2024leveraging, fornasier2024consensus, gerber2025mean}.

Global optimization problems can be classified into two main categories depending on the countability of $\Omega$: if $\Omega$ is finite or countably infinite, the problem is called {\it discrete global optimization}; if $\Omega$ is uncountable, the problem is called {\it continuous global optimization}. The strategies used in these two categories are very different, and in this paper we focus on continuous global optimization. In this setting, the function $f$ is usually assumed to be continuous, and the search space $\Omega$ is often assumed to be compact, which ensures that $f$ has a global maximum and minimum within $\Omega$. Using equivalently either the $\arg\max$ or the $\arg\min$ operator, the problem of continuous global optimization can be formally expressed as:
$$
  x^\star \in \argmax{x \in \Omega} f(x) \ = \ \argmin{x \in \Omega} (-f(x)).
$$

In the special case when $f$ is a convex function, the problem of global optimization is equivalent to that of local optimization, which explains why global optimization is often called {\it non-convex optimization}. The global optimization strategies are different from those used to solve local optimization problems, and the theoretical analysis of the former algorithms can be more challenging, as the class of functions considered is much larger and the algorithms often explore stochastically the search space \cite{Hansen1996, Pinnau2017, Frazier2018, Contal2013, Serre2025}.


\inlinetitle{Existing global optimization libraries and limitations}{.}
The landscape of global optimization software is diverse, with many libraries addressing specific problems or focusing on algorithm families. We review the most prominent existing solutions:

\noindent\emph{\mySqBullet\ NLopt \cite{NLopt}}: Provides a unified interface to classical algorithms (DIRECT, CRS, evolutionary strategies) but lacks structural pattern exploitation and recent advances.

\noindent\emph{\mySqBullet\ Scipy.optimize \cite{Virtanen2020}}: Broad optimization library where global optimization is secondary to local methods.

\noindent\emph{\mySqBullet\ Bayesian Optimization \cite{Nogueira2014}}: Specialized for expensive black-box functions, not general-purpose.

\noindent\emph{\mySqBullet\ Optuna \cite{Akiba2019} and Ray Tune \cite{Liaw2018}}: Both specialized for hyperparameter tuning in machine learning rather than general continuous global optimization.

While existing libraries serve their respective domains well, each specializes in particular optimization scenarios, and therefore leave a gap in the available open-source solutions:

\begin{enumerate}[leftmargin=1.5em]
  \item \textbf{Specialization vs. generality:} Most libraries focus on a specific problem type (\eg hyperparameter tuning, expensive black-box functions) or a specific algorithm family (\eg evolutionary algorithms, Bayesian optimization). A unified framework for general continuous global optimization is lacking.
  
  \item \textbf{Limited code modularity:} Algorithms are typically implemented independently without exploiting common patterns, such as sampling strategies, population evolution, decision making loops. This leads to code duplication and makes it difficult to add new algorithms or features that span multiple families.
  
  \item \textbf{Missing recent advances:} Many libraries focus on classical algorithms. Recent advances in stochastic differential equation (SDE)-driven methods, particle-based algorithms with common noise, and modern decision-based approaches are rarely integrated. For research purposes, comparing newly proposed methods against these recent advances is essential to validate their effectiveness, yet existing libraries do not provide easy access to such state-of-the-art methods.
\end{enumerate}

\inlinetitle{Contribution of the \globe library}{.}
\globe directly addresses these challenges through a modular, general-purpose framework for continuous global optimization that: \textbf{i)} consolidates classical and state-of-the-art algorithms, enabling researchers to benchmark new methods against up-to-date baselines; \textbf{ii)} provides integrated benchmarking tools for direct and reproducible algorithm comparison; and \textbf{iii)} factors common patterns at the family level to facilitate the creation of new algorithms and plugins. This modular architecture is grounded in recent advances in the mathematical formalization of global optimization \cite{Powell2019, Serre2026Unifying, Serre2026}, which revealed structural commonalities across algorithms and enabled the design of broadly applicable, formally justified features. By reflecting these mathematical structures, \globe achieves both extensibility and research-friendliness. The C++ backend relies on the Eigen library to provide efficient vectorized linear algebra operations, ensuring high performance across all algorithms. \globe currently includes $14$ optimizers, $19$ analytical benchmarks with a random function generator, and an integrated toolkit for algorithm comparison.

The documentation of \globe is available online\footurl{https://gaetanserre.fr/GLOBe} and the code is open-source under the GPL-3.0 license\footurl{https://github.com/gaetanserre/GLOBe}.

\section{The \globe library}

\globe currently implements $14$ global optimization algorithms, as summarized in \cref{tbl:algorithms}. 

\subsection{Modular Architecture}

\begin{figure}[tb]
  \includegraphics[width=\linewidth]{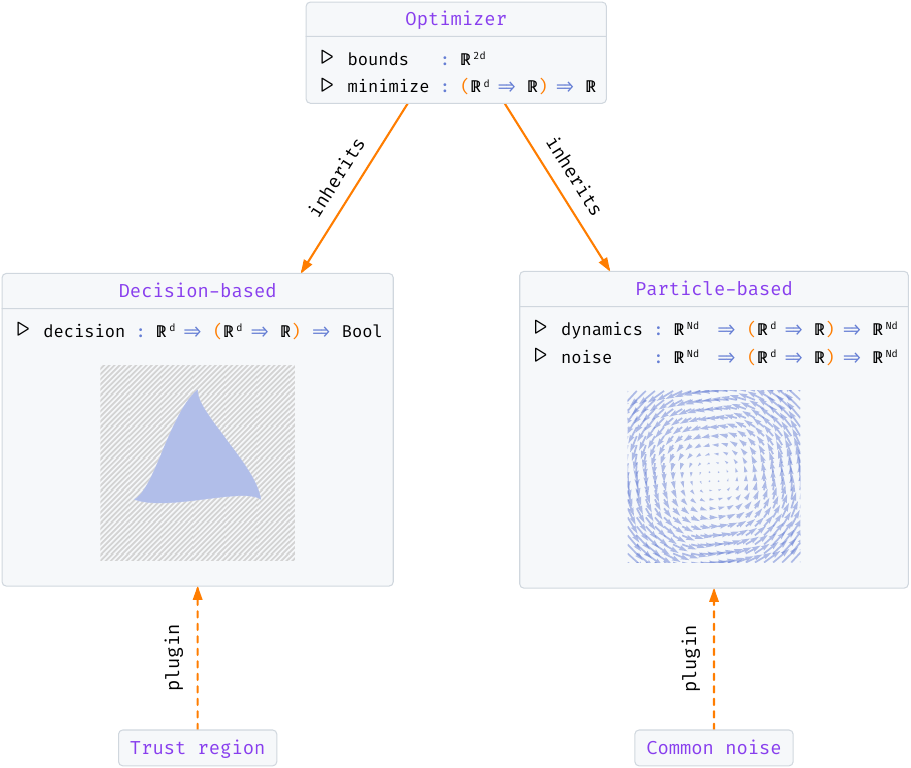}
  \caption{\globe's architecture showing the class hierarchy and modular decomposition. The base \texttt{Optimizer} interface is extended by concrete implementations organized in two principal families: \texttt{Decision} (for decision-based methods) and \texttt{Particles\_Optimizer} (for particle-based methods). Orthogonal modules, such as trust regions and common noise, can be composed with these base classes to create algorithm variants.}
  \label{fig:arch}
\end{figure}

The core design philosophy of \globe is to achieve \emph{composition through inheritance and modularity}. Rather than treating all algorithms as independent implementations, \globe identifies common patterns across algorithm families and factors them into reusable components. This family-level factorization is especially important for plugins: modules are attached to family base classes, so a single plugin implementation is valid for all algorithms in that family. This approach is illustrated in \cref{fig:arch}.

\begin{table}[t]
\centering
\caption{Optimization algorithms implemented in \globe.}
\begin{adjustbox}{width=\linewidth}
\begin{tabular}{lll}
\toprule
\textbf{Algorithm} & \textbf{Family} & \textbf{Reference} \\
\toprule
AdaLIPO+ & Decision-based & \citet{Serre2024} \\
ECP (Every Call is Precious) & Decision-based & \citet{Fourati2025} \\
AdaRankOpt & Decision-based & \citet{Malherbe2016} \\
\hline
PSO (Particle Swarm) & Particle-based & \citet{Kennedy1995} \\
CBO (Consensus-Based Opt.) & Particle-based & \citet{Pinnau2017} \\
Langevin Dynamics & Particle-based & \citet{Prigogine1957} \\
SBS (Stein Boltzmann Sampling) & Particle-based & \citet{Serre2025} \\
MSGD (Multi-Start SGD) & Particle-based & -- \\
\hline
DIRECT & Misc. & \citet{Jones2001} \\
CRS (Controlled Random Search) & Misc. & \citet{Price1983} \\
CMA-ES & Misc. & \citet{Hansen1996} \\
MLSL (Multi-Level Single-Linkage) & Misc. & \citet{Rinnooy1985} \\
Gradient Descent & Misc. & -- \\
Pure Random Search & Misc. & -- \\
\bottomrule
\end{tabular}
\end{adjustbox}
\label{tbl:algorithms}
\end{table}

This modular architecture goes beyond engineering convenience, as it is deeply rooted in recent advances in the mathematical formalization of global optimization algorithms. Works such as \cite{Powell2019, Serre2026Unifying, Serre2026} have revealed fundamental structural patterns common to multiple algorithm families, formally characterizing the mathematical invariants and common building blocks that transcend individual algorithms. Notably, \citet{Serre2026Unifying} develops frameworks enabling the application of formal methods to global optimization, providing a principled foundation for algorithm analysis and composition. By codifying these formal structures, \globe's modular design reflects the mathematical reality of global optimization, enabling a principled approach to algorithm development and composition. The adopted perspective aligns with the broader trend toward formal methods in machine learning and optimization research, where structured mathematical frameworks provide both theoretical grounding and practical extensibility.

\subsection{Base Optimizer Interface}

\begin{figure}[tb]
  \centering
  \includegraphics{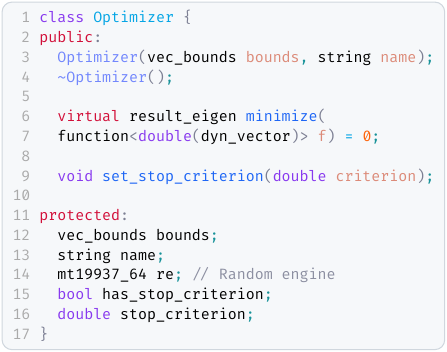}
	 \caption{Abstract base \texttt{Optimizer} class interface showing the minimal API that all algorithms must implement.}
  \label{fig:interface}
\end{figure}

All optimization algorithms in \globe inherit from a common abstract base class that defines the minimal interface required, as shown in \cref{fig:interface}.
This interface ensures that all algorithms can be used interchangeably. The base class handles:
\begin{itemize}[leftmargin=1.5em]
  \item \textbf{Bounds management}: The search space definition.
  \item \textbf{Random number generation}: Seeded random number generation engine for reproducibility.
  \item \textbf{Stopping criteria}: Pluggable stopping conditions (target values).
  \item \textbf{Python interoperability}: Provides the C++/Python interface layer. As a result, this interoperability burden is already handled when implementing new algorithms in C++.
\end{itemize}

\subsection{Decision-Based Optimizers}

Decision-based algorithms make progress by maintaining a history of function evaluations, and according to this history they decide which region of the search space to sample from next. Examples of this class include AdaLIPO+ (Lipschitz-based adaptive partitioning) \cite{Serre2024}, ECP (Every Call is Precious) \cite{Fourati2025}, and AdaRankOpt (ranking-based decisions) \cite{Malherbe2016}.

The common pattern across these algorithms is:
\begin{enumerate}
  \item \textit{Initialize}: Sample an initial point $x_0 \in \Omega$ uniformly at random.
  \item \textit{Iterate}: At iteration $t$, using the history $\{(x_i, f(x_i))\}_{i=0}^{t-1}$, compute a region $R_t \subseteq \Omega$.
  \item \textit{Sample}: Draw a new candidate point $x_t$ at random over $R_t$.
  \item \textit{Evaluate}: Compute $f(x_t)$ and append it to history.
\end{enumerate}

This pattern is defined by the \texttt{Decision} base class that inherits from \texttt{Optimizer}. The \texttt{Decision} class handles the evaluation loop and history management, while specific algorithms (AdaLIPO+, ECP, etc.) only need to implement their decision-making logic via callback functions. This factorization avoids duplicating the iteration framework and sampling process across multiple decision-based methods, and allows family-level plugins (\eg trust-region refinement) to be implemented once and reused by all decision-based algorithms.

\subsection{Particle-Based Optimizers}

Particle-based methods maintain a population of $N$ particles that evolve at each iteration according to stochastic dynamics. Many such algorithms can be modeled as discretizations of McKean-Vlasov stochastic differential equations:
\begin{equation}\label{eq:mckean-vlasov}
\dif X^{(i)}_t = b(X^{(i)}_t, \hat{\mu}_t) \dif t + \sigma(X^{(i)}_t, \hat{\mu}_t) \dif B^{(i)}_t,
\end{equation}
where $b$ is the drift (deterministic evolution), $\sigma\,\dif B^{(i)}_t$ is the diffusion (stochastic noise), and $\hat{\mu}_t$ is the empirical law of $X^{(i)}_t$. Examples of particle-based algorithms in \globe include Particle Swarm Optimization (PSO) \cite{Kennedy1995}, Consensus-Based Optimization (CBO) \cite{Pinnau2017}, Langevin dynamics \cite{Prigogine1957}, Stein Boltzmann Sampling (SBS) \cite{Serre2025}, and Multi-Start Gradient Descent (MSGD).

The common pattern is:
\begin{enumerate}[leftmargin=1.5em]
  \item \textit{Initialize}: Draw $N$ particles $X^{(1)}_0, \ldots, X^{(N)}_0$ uniformly from $\Omega$.
  \item \textit{Iterate}: At each iteration, compute drift and noise terms, then update all particles.
  \item \textit{Discretize}: Use an Euler-Maruyama scheme to discretize the SDE and update particles.
\end{enumerate}

This pattern is defined by the \texttt{Particles\_Optimizer} base class that inherits from \texttt{Optimizer}. Concrete algorithms (PSO, CBO, SBS, \ldots) only need to implement the drift and noise computation for their specific dynamics. The base class handles the Euler-Maruyama discretization loop, particle updates, boundary handling, and adaptive scheduling; it also provides a natural integration point for family-level plugins (\eg common noise), so a plugin implementation can be applied consistently across all particle-based algorithms.

\subsection{Modular Plugins}

Beyond the base hierarchies, \globe provides orthogonal modules that enhance algorithms. \\

\inlinetitle{Trust Region Module}{.~}
Decision-based algorithms can be augmented with local refinement using trust-region methods. Specifically, when a promising region is identified, the algorithm can temporarily switch to a local optimizer (BOBYQA \cite{Powell2009}) to refine the solution before resuming global exploration. This is implemented as an optional family-level module on top of the \texttt{Decision} base class, making the plugin valid for all decision-based optimizers without modifying their core logic.

\inlinetitle{Common Noise Module}{.~}
Particle-based methods based on McKean-Vlasov dynamics face a fundamental challenge: when particles evolve under independent noise alone, exploration efficiency suffers due to insufficient diversity in particle trajectories \cite{Tugaut2014}. The recent work in \cite{Serre2026}, addresses this by introducing common noise schemes that are shared across all particles, enhancing global exploration:
\begin{itemize}[leftmargin=1.5em]
  \item \emph{SMD (Stochastic Moment Dynamics):} Injects noise into aggregate particle statistics, specifically the population mean and variance. This approach perturbs the collective evolution, improving coverage of the search space.
  
  \item \emph{GCN (Geometrical Common Noise):} Constructs a canonical Gaussian common noise using the geometric structure of the underlying process.
\end{itemize}

\globe implements these noise schemes as modular components integrated at the \texttt{Particles\_Optimizer} family level, so the same plugin is valid for any particle-based optimizer. This allows users to enhance exploration capabilities without modifying algorithm-specific core dynamics.

\inlinetitle{Particle Filtering}{.~}
\globe provides a flexible filtering interface at the family level, enabling algorithm-specific strategies for all particle-based optimizers. Currently, it implements quantile-based filtering \cite{Serre2025} that removes particles with small displacement and poor function values, improving exploration efficiency. New strategies can be added without modifying core algorithm dynamics such as Sequential Importance Resampling \cite{Kitagawa1993}.

\section{Benchmarking Features}

\subsection{Benchmarking Functions}

\globe includes a wide range of synthetic benchmark functions, as well as a random generator of functions, which can be used to evaluate the performance of optimization algorithms. It also implements a toolkit to perform large-scale experiments and analyze results. All benchmarks can run in any arbitrary number of dimensions.

\inlinetitle{Benchmark Interface}{.~}%
All benchmark functions inherit from a common \texttt{Benchmark} class that provides a standardized interface. Each benchmark encapsulates the function expression, tracks evaluation counts, provides gradient estimates via finite differences, and supports visualization. This unified API enables consistent comparison across all benchmark types.

\inlinetitle{Analytical Benchmark Functions}{.~}%
\globe includes a comprehensive suite of $19$ analytical benchmark functions organized by modality, see \cref{tbl:benchmarks}.
Most functions include known global minima, allowing precise assessment of algorithm performance.

\begin{table}[t]
\centering
\caption{Analytical benchmark functions in \globe.}
\begin{adjustbox}{width=\linewidth}
\begin{tabular}{lll}
\toprule
\textbf{Category} & \textbf{Functions} & \textbf{Characteristics} \\
\toprule
Unimodal & 
\begin{tabular}[c]{@{}l@{}} Sphere, Ellipsoid, Bentcigar, \\ Hyperellipsoid, Zakharov, \\ Trid, Sumpow, Rosenbrock, \\ Dixonprice \end{tabular}
&
\begin{tabular}[c]{@{}l@{}} Single optimum, smooth \\ landscapes, test convergence \\  speed \end{tabular} \\
\hline
Multimodal &
\begin{tabular}[c]{@{}l@{}} Rastrigin, Ackley, Schwefel, \\ Levy, Michalewicz, Langermann, \\ Deb, Griewank, \\ Styblinski-Tang \end{tabular}
&
\begin{tabular}[c]{@{}l@{}} Numerous local minima, \\  test global exploration \end{tabular} \\
\bottomrule
\end{tabular}
\end{adjustbox}
\label{tbl:benchmarks}
\end{table}

\inlinetitle{Random Function Generator}{.~}%
\globe integrates PyGKLS, a Python wrapper around the GKLS function generator~\cite{Gaviano2003}. GKLS generates random smooth multidimensional functions with:
\begin{itemize}[leftmargin=1.5em]
  \item known global and local minima locations;
  \item controllable multi-modality degree;
  \item controllable smoothness (from $C^2$ to $C^0$);
  \item reproducibility via fixed seeds.
\end{itemize}

The GKLS generator provides a rigorous alternative to fixed analytical benchmarks, enabling statistical analysis across random function instances.

\subsection{Benchmarking Engine}

\globe provides an integrated benchmarking framework via its main \texttt{GLOBe} class, which orchestrates comprehensive experimental comparisons. The class takes as input three main components: a list of optimizers (with optional hyperparameter configurations), a list of benchmark functions, and a list of metrics. This modular design enables systematic analysis across multiple algorithms, problems, and performance criteria.

The typical workflow is as follows:
\begin{enumerate}[leftmargin=1.5em]
  \item \emph{Configuration}: Users specify optimizer instances (or names with hyperparameters), benchmark functions, and evaluation metrics.
  \item \emph{Multiple runs}: The framework runs each optimizer-benchmark combination $n$ times to account for stochasticity, and collects the solution values from each run.
  \item \emph{Metric computation}: For each optimizer-benchmark pair, the framework computes all requested metrics across the collected solutions. This enables richer analysis than simple optimality gaps.
  \item \emph{Aggregation and reporting}: Results are aggregated and displayed as comparative tables (console or \LaTeX{} format), allowing side-by-side algorithm comparison.
\end{enumerate}

\inlinetitle{Extensible Metrics}{.~}
\globe provides a flexible framework where each metric is a callable that processes the solutions from multiple runs. Built-in metrics include:
\begin{itemize}[leftmargin=1.5em]
  \item \emph{Approximation quality}: Mean and standard deviation of solution values across runs, providing robustness assessment.
  \item \emph{$f$-target value}: A problem-adaptive threshold computed as $f_{\text{target}} = f_{\min} + (f_{\text{mean}} - f_{\min})(1-p)$, where $f_{\min}$ is the global minimum, $f_{\text{mean}}$ is the mean value over a large random sample from $\Omega$, and $p$ is a configurable proportion (default $p=0.99$). This represents the quantile of the random distribution of function values.
  \item \emph{Success rate (proportion)}: Percentage of runs where the algorithm achieves a solution value below the $f$-target threshold, effectively measuring how often the optimizer succeeds in reaching a specified percentile performance level.
\end{itemize}

The metrics interface is designed to be extensible: users can implement custom metrics by inheriting from the base \texttt{Metric} class and by defining their evaluation logic. This enables domain-specific performance analysis (\eg computational budget constraints, constraint satisfaction for engineering problems).

\section{Conclusion}

We presented \globe, a modular general-purpose framework addressing fragmentation in global optimization software. By consolidating classical and recent state-of-the-art algorithms into a unified platform, \globe provides researchers with immediate access to diverse methods and enables direct comparisons. The modular architecture is grounded in recent advances in the mathematical formalization of global optimization \cite{Powell2019, Serre2026Unifying, Serre2026}, which reveal structural patterns across algorithm families. By codifying these patterns as reusable components, \globe reduces code duplication and facilitates algorithm development. Family-level plugins such as trust-region refinement, common noise enhancement, and particle filtering are implemented once and automatically available to all algorithms in their families. The integrated benchmarking toolkit with its extensible metrics framework enables rigorous reproducible comparisons. By bridging mathematical formalization and software design, \globe serves both research and practice, providing a practical foundation for continuous global optimization.

\section{Acknowledgments}
This work was supported by the Industrial Data Analytics and Machine Learning Chair hosted at ENS Paris-Saclay.

\bibliographystyle{ACM-Reference-Format}
\bibliography{refs}


\begin{thebibliography}{37}


\ifx \showCODEN    \undefined \def \showCODEN     #1{\unskip}     \fi
\ifx \showISBNx    \undefined \def \showISBNx     #1{\unskip}     \fi
\ifx \showISBNxiii \undefined \def \showISBNxiii  #1{\unskip}     \fi
\ifx \showISSN     \undefined \def \showISSN      #1{\unskip}     \fi
\ifx \showLCCN     \undefined \def \showLCCN      #1{\unskip}     \fi
\ifx \shownote     \undefined \def \shownote      #1{#1}          \fi
\ifx \showarticletitle \undefined \def \showarticletitle #1{#1}   \fi
\ifx \showURL      \undefined \def \showURL       {\relax}        \fi
\providecommand\bibfield[2]{#2}
\providecommand\bibinfo[2]{#2}
\providecommand\natexlab[1]{#1}
\providecommand\showeprint[2][]{arXiv:#2}

\bibitem[Agrawal(2021)]%
        {Agrawal2021}
\bibfield{author}{\bibinfo{person}{Tanay Agrawal}.}
  \bibinfo{year}{2021}\natexlab{}.
\newblock \showarticletitle{Hyperparameter optimization in machine learning:
  make your machine learning and deep learning models more efficient}.
\newblock \bibinfo{journal}{\emph{New York, NY}} (\bibinfo{year}{2021}),
  \bibinfo{pages}{109--129}.
\newblock


\bibitem[Akiba et~al\mbox{.}(2019)]%
        {Akiba2019}
\bibfield{author}{\bibinfo{person}{Takuya Akiba}, \bibinfo{person}{Shotaro
  Sano}, \bibinfo{person}{Toshihiko Yanase}, \bibinfo{person}{Takeru Ohta},
  {and} \bibinfo{person}{Masanori Koyama}.} \bibinfo{year}{2019}\natexlab{}.
\newblock \showarticletitle{Optuna: A next-generation hyperparameter
  optimization framework}. In \bibinfo{booktitle}{\emph{Proceedings of the 25th
  ACM SIGKDD \International conference on knowledge discovery \& data mining}}.
\newblock


\bibitem[Awasthi et~al\mbox{.}(2024)]%
        {Awasthi2024}
\bibfield{author}{\bibinfo{person}{Pranjal Awasthi}, \bibinfo{person}{Anqi
  Mao}, \bibinfo{person}{Mehryar Mohri}, {and} \bibinfo{person}{Yutao Zhong}.}
  \bibinfo{year}{2024}\natexlab{}.
\newblock \showarticletitle{{DC}-programming for neural network optimizations}.
\newblock \bibinfo{journal}{\emph{\Journal of Global Optimization}}
  (\bibinfo{year}{2024}).
\newblock


\bibitem[Bungert et~al\mbox{.}(2024)]%
        {Bungert2024}
\bibfield{author}{\bibinfo{person}{Leon Bungert}, \bibinfo{person}{Tim Roith},
  {and} \bibinfo{person}{Philipp Wacker}.} \bibinfo{year}{2024}\natexlab{}.
\newblock \showarticletitle{Polarized consensus-based dynamics for optimization
  and sampling}.
\newblock \bibinfo{journal}{\emph{Mathematical Programming}}
  \bibinfo{volume}{211} (\bibinfo{year}{2024}).
\newblock


\bibitem[Contal et~al\mbox{.}(2013)]%
        {Contal2013}
\bibfield{author}{\bibinfo{person}{Emile Contal}, \bibinfo{person}{David
  Buffoni}, \bibinfo{person}{Alexandre Robicquet}, {and}
  \bibinfo{person}{Nicolas Vayatis}.} \bibinfo{year}{2013}\natexlab{}.
\newblock \showarticletitle{Parallel {G}aussian process optimization with upper
  confidence bound and pure exploration}. In \bibinfo{booktitle}{\emph{Joint
  European \Conference on Machine Learning and Knowledge Discovery in
  Databases}}.
\newblock


\bibitem[El~Yadari et~al\mbox{.}(2024)]%
        {el2024taxonomy}
\bibfield{author}{\bibinfo{person}{Meryeme El~Yadari}, \bibinfo{person}{Saloua
  El~Motaki}, \bibinfo{person}{Ali Yahyaouy}, \bibinfo{person}{Philippe
  Makany}, \bibinfo{person}{Khalid El~Fazazy}, \bibinfo{person}{Hamid Gualous},
  {and} \bibinfo{person}{St{\'e}phane Le~Masson}.}
  \bibinfo{year}{2024}\natexlab{}.
\newblock \showarticletitle{Taxonomy of optimization algorithms combined with
  {CNN} for optimal placement of virtual machines within physical machines in
  data centers}.
\newblock \bibinfo{journal}{\emph{Energy Informatics}} (\bibinfo{year}{2024}).
\newblock


\bibitem[Fornasier et~al\mbox{.}(2024)]%
        {fornasier2024consensus}
\bibfield{author}{\bibinfo{person}{Massimo Fornasier}, \bibinfo{person}{Timo
  Klock}, {and} \bibinfo{person}{Konstantin Riedl}.}
  \bibinfo{year}{2024}\natexlab{}.
\newblock \showarticletitle{Consensus-based optimization methods converge
  globally}.
\newblock \bibinfo{journal}{\emph{SIAM \Journal on Optimization}}
  (\bibinfo{year}{2024}).
\newblock


\bibitem[Fourati et~al\mbox{.}(2025)]%
        {Fourati2025}
\bibfield{author}{\bibinfo{person}{Fares Fourati}, \bibinfo{person}{Salma
  Kharrat}, \bibinfo{person}{Vaneet Aggarwal}, {and}
  \bibinfo{person}{Mohamed-Slim Alouini}.} \bibinfo{year}{2025}\natexlab{}.
\newblock \showarticletitle{{Every Call is Precious}: Global Optimization of
  Black-Box Functions with Unknown Lipschitz Constants}.
\newblock \bibinfo{journal}{\emph{\International \Conference on Artificial
  Intelligence and Statistics}}.
\newblock


\bibitem[Franceschi et~al\mbox{.}(2024)]%
        {Franceschi2024}
\bibfield{author}{\bibinfo{person}{Luca Franceschi}, \bibinfo{person}{Michele
  Donini}, \bibinfo{person}{Valerio Perrone}, \bibinfo{person}{Aaron Klein},
  \bibinfo{person}{C{\'e}dric Archambeau}, \bibinfo{person}{Matthias Seeger},
  \bibinfo{person}{Massimiliano Pontil}, {and} \bibinfo{person}{Paolo
  Frasconi}.} \bibinfo{year}{2024}\natexlab{}.
\newblock \showarticletitle{Hyperparameter optimization in machine learning}.
\newblock \bibinfo{journal}{\emph{Preprint arXiv:2410.22854}}
  (\bibinfo{year}{2024}).
\newblock


\bibitem[Frazier(2018)]%
        {Frazier2018}
\bibfield{author}{\bibinfo{person}{Peter~I. Frazier}.}
  \bibinfo{year}{2018}\natexlab{}.
\newblock \bibinfo{title}{A Tutorial on {B}ayesian Optimization}.
\newblock


\bibitem[Gaviano et~al\mbox{.}(2003)]%
        {Gaviano2003}
\bibfield{author}{\bibinfo{person}{Marco Gaviano}, \bibinfo{person}{Dmitri~E.
  Kvasov}, \bibinfo{person}{Daniela Lera}, {and} \bibinfo{person}{Yaroslav~D.
  Sergeyev}.} \bibinfo{year}{2003}\natexlab{}.
\newblock \showarticletitle{Algorithm 829: Software for generation of classes
  of test functions with known local and global minima for global
  optimization}.
\newblock \bibinfo{journal}{\emph{ACM Trans. Math. Software}}
  (\bibinfo{year}{2003}).
\newblock


\bibitem[Gerber et~al\mbox{.}(2025)]%
        {gerber2025mean}
\bibfield{author}{\bibinfo{person}{Nicolai~Jurek Gerber},
  \bibinfo{person}{Franca Hoffmann}, {and} \bibinfo{person}{Urbain Vaes}.}
  \bibinfo{year}{2025}\natexlab{}.
\newblock \showarticletitle{Mean-field limits for consensus-based optimization
  and sampling}.
\newblock \bibinfo{journal}{\emph{ESAIM: Control, Optimisation and Calculus of
  Variations}}  \bibinfo{volume}{31} (\bibinfo{year}{2025}),
  \bibinfo{pages}{74}.
\newblock


\bibitem[Hansen and Ostermeier(1996)]%
        {Hansen1996}
\bibfield{author}{\bibinfo{person}{Nikolaus Hansen} {and}
  \bibinfo{person}{Andreas Ostermeier}.} \bibinfo{year}{1996}\natexlab{}.
\newblock \showarticletitle{Adapting arbitrary normal mutation distributions in
  evolution strategies: The covariance matrix adaptation}. In
  \bibinfo{booktitle}{\emph{IEEE \International \Conference on Evolutionary
  Computation}}.
\newblock


\bibitem[Houssein et~al\mbox{.}(2024)]%
        {Houssein2024}
\bibfield{author}{\bibinfo{person}{Essam~H. Houssein},
  \bibinfo{person}{Mahmoud~Khalaf Saeed}, \bibinfo{person}{Gang Hu}, {and}
  \bibinfo{person}{Mustafa~M. Al-Sayed}.} \bibinfo{year}{2024}\natexlab{}.
\newblock \showarticletitle{Metaheuristics for Solving Global and Engineering
  Optimization Problems: Review, Applications, Open Issues and Challenges}.
\newblock \bibinfo{journal}{\emph{Archives of Computational Methods in
  Engineering}} (\bibinfo{year}{2024}).
\newblock


\bibitem[Johnson(2007)]%
        {NLopt}
\bibfield{author}{\bibinfo{person}{Steven~G. Johnson}.}
  \bibinfo{year}{2007}\natexlab{}.
\newblock \bibinfo{title}{The {NLopt} nonlinear-optimization package}.
\newblock \bibinfo{howpublished}{\url{https://github.com/stevengj/nlopt}}.
\newblock


\bibitem[Jones(2001)]%
        {Jones2001}
\bibfield{author}{\bibinfo{person}{Donald~R. Jones}.}
  \bibinfo{year}{2001}\natexlab{}.
\newblock \bibinfo{booktitle}{\emph{Direct Global Optimization Algorithm}}.
\newblock


\bibitem[Kennedy and Eberhart(1995)]%
        {Kennedy1995}
\bibfield{author}{\bibinfo{person}{James Kennedy} {and}
  \bibinfo{person}{Russell~C. Eberhart}.} \bibinfo{year}{1995}\natexlab{}.
\newblock \showarticletitle{Particle swarm optimization}. In
  \bibinfo{booktitle}{\emph{\International \Conference on Neural Networks}}.
\newblock


\bibitem[Kirkpatrick et~al\mbox{.}(1983)]%
        {Kirkpatrick1983}
\bibfield{author}{\bibinfo{person}{Scott Kirkpatrick},
  \bibinfo{person}{C~Daniel Gelatt~Jr}, {and} \bibinfo{person}{Mario~P
  Vecchi}.} \bibinfo{year}{1983}\natexlab{}.
\newblock \showarticletitle{Optimization by Simulated Annealing}.
\newblock \bibinfo{journal}{\emph{Science}} (\bibinfo{year}{1983}).
\newblock


\bibitem[Kitagawa(1993)]%
        {Kitagawa1993}
\bibfield{author}{\bibinfo{person}{Genshiro Kitagawa}.}
  \bibinfo{year}{1993}\natexlab{}.
\newblock \showarticletitle{A Monte Carlo filtering and smoothing method for
  non-Gaussian nonlinear state space models}. In
  \bibinfo{booktitle}{\emph{Proceedings of the 2nd US-Japan joint seminar on
  statistical time series analysis}}.
\newblock


\bibitem[Lee et~al\mbox{.}(2017)]%
        {lee2017finding}
\bibfield{author}{\bibinfo{person}{Juyong Lee}, \bibinfo{person}{In-Ho Lee},
  \bibinfo{person}{InSuk Joung}, \bibinfo{person}{Jooyoung Lee}, {and}
  \bibinfo{person}{Bernard~R Brooks}.} \bibinfo{year}{2017}\natexlab{}.
\newblock \showarticletitle{Finding multiple reaction pathways via global
  optimization of action}.
\newblock \bibinfo{journal}{\emph{Nature Communications}}
  (\bibinfo{year}{2017}).
\newblock


\bibitem[Liaw et~al\mbox{.}(2018)]%
        {Liaw2018}
\bibfield{author}{\bibinfo{person}{Richard Liaw}, \bibinfo{person}{Eric Liang},
  \bibinfo{person}{Robert Nishihara}, \bibinfo{person}{Philipp Moritz},
  \bibinfo{person}{Joseph~E. Gonzalez}, {and} \bibinfo{person}{Ion Stoica}.}
  \bibinfo{year}{2018}\natexlab{}.
\newblock \bibinfo{title}{Tune: A Research Platform for Distributed Model
  Selection and Training}.
\newblock


\bibitem[Malherbe et~al\mbox{.}(2016)]%
        {Malherbe2016}
\bibfield{author}{\bibinfo{person}{Cedric Malherbe}, \bibinfo{person}{Emile
  Contal}, {and} \bibinfo{person}{Nicolas Vayatis}.}
  \bibinfo{year}{2016}\natexlab{}.
\newblock \showarticletitle{A ranking approach to global optimization}. In
  \bibinfo{booktitle}{\emph{\International \Conference on Machine Learning}}.
\newblock


\bibitem[Nogueira(14  )]%
        {Nogueira2014}
\bibfield{author}{\bibinfo{person}{Fernando Nogueira}.}
  \bibinfo{year}{2014--}\natexlab{}.
\newblock \bibinfo{title}{{Bayesian Optimization}: Open source constrained
  global optimization tool for {Python}}.
\newblock
\urldef\tempurl%
\url{https://github.com/bayesian-optimization/BayesianOptimization}
\showURL{%
\tempurl}


\bibitem[Pinnau et~al\mbox{.}(2017)]%
        {Pinnau2017}
\bibfield{author}{\bibinfo{person}{Ren{\'e} Pinnau}, \bibinfo{person}{Claudia
  Totzeck}, \bibinfo{person}{Oliver Tse}, {and} \bibinfo{person}{Stephan
  Martin}.} \bibinfo{year}{2017}\natexlab{}.
\newblock \showarticletitle{A consensus-based model for global optimization and
  its mean-field limit}.
\newblock \bibinfo{journal}{\emph{Mathematical Models and Methods in Applied
  Sciences}} (\bibinfo{year}{2017}).
\newblock


\bibitem[Pintér(1996)]%
        {Pinter1996}
\bibfield{author}{\bibinfo{person}{János~D. Pintér}.}
  \bibinfo{year}{1996}\natexlab{}.
\newblock \showarticletitle{Global Optimization in Action}.
\newblock \bibinfo{journal}{\emph{Nonconvex Optimization and Its Applications}}
  (\bibinfo{year}{1996}).
\newblock


\bibitem[Powell et~al\mbox{.}(2009)]%
        {Powell2009}
\bibfield{author}{\bibinfo{person}{Michael~JD Powell} {et~al\mbox{.}}}
  \bibinfo{year}{2009}\natexlab{}.
\newblock \showarticletitle{The {BOBYQA} algorithm for bound constrained
  optimization without derivatives}.
\newblock \bibinfo{journal}{\emph{Cambridge NA Report NA2009/06, University of
  Cambridge, Cambridge}} (\bibinfo{year}{2009}).
\newblock


\bibitem[Powell(2019)]%
        {Powell2019}
\bibfield{author}{\bibinfo{person}{Warren~B Powell}.}
  \bibinfo{year}{2019}\natexlab{}.
\newblock \showarticletitle{A unified framework for stochastic optimization}.
\newblock \bibinfo{journal}{\emph{European Journal of Operational Research}}
  (\bibinfo{year}{2019}).
\newblock


\bibitem[Price(1983)]%
        {Price1983}
\bibfield{author}{\bibinfo{person}{W.~L. Price}.}
  \bibinfo{year}{1983}\natexlab{}.
\newblock \showarticletitle{Global optimization by controlled random search}.
\newblock \bibinfo{journal}{\emph{\Journal of Optimization Theory and
  Applications}} (\bibinfo{year}{1983}).
\newblock


\bibitem[Prigogine and Balescu(1957)]%
        {Prigogine1957}
\bibfield{author}{\bibinfo{person}{Ilya Prigogine} {and} \bibinfo{person}{Radu
  Balescu}.} \bibinfo{year}{1957}\natexlab{}.
\newblock \showarticletitle{Sur la theorie moleculaire du mouvement brownien}.
\newblock \bibinfo{journal}{\emph{Physica}} (\bibinfo{year}{1957}).
\newblock


\bibitem[Riedl(2024)]%
        {riedl2024leveraging}
\bibfield{author}{\bibinfo{person}{Konstantin Riedl}.}
  \bibinfo{year}{2024}\natexlab{}.
\newblock \showarticletitle{Leveraging memory effects and gradient information
  in consensus-based optimisation: On global convergence in mean-field law}.
\newblock \bibinfo{journal}{\emph{European \Journal of Applied Mathematics}}
  (\bibinfo{year}{2024}).
\newblock


\bibitem[Rinnooy~Kan et~al\mbox{.}(1985)]%
        {Rinnooy1985}
\bibfield{author}{\bibinfo{person}{A.~H.~G. Rinnooy~Kan},
  \bibinfo{person}{G.~T. Timmer}, \bibinfo{person}{A.~H.~G. Rinnooy~Kan}, {and}
  \bibinfo{person}{G.~T. Timmer}.} \bibinfo{year}{1985}\natexlab{}.
\newblock \showarticletitle{The Multi Level Single Linkage Method for
  Unconstrained and Constrained Global Optimization}.
\newblock  (\bibinfo{year}{1985}).
\newblock


\bibitem[Serr\'{e} et~al\mbox{.}(2024)]%
        {Serre2024}
\bibfield{author}{\bibinfo{person}{Ga\"{e}tan Serr\'{e}},
  \bibinfo{person}{Perceval Beja-Battais}, \bibinfo{person}{Sophia Chirrane},
  \bibinfo{person}{Argyris Kalogeratos}, {and} \bibinfo{person}{Nicolas
  Vayatis}.} \bibinfo{year}{2024}\natexlab{}.
\newblock \showarticletitle{{LIPO}+: Frugal Global Optimization for Lipschitz
  Functions}. In \bibinfo{booktitle}{\emph{Hellenic \Conference on Artificial
  Intelligence}}.
\newblock


\bibitem[Serré et~al\mbox{.}(2026a)]%
        {Serre2026}
\bibfield{author}{\bibinfo{person}{Gaëtan Serré}, \bibinfo{person}{Pierre
  Germain}, \bibinfo{person}{Samuel Gruffaz}, {and} \bibinfo{person}{Argyris
  Kalogeratos}.} \bibinfo{year}{2026}\natexlab{a}.
\newblock \showarticletitle{Enhancing Exploration in Global Optimization by
  Noise Injection in the Probability Measures Space}.
\newblock \bibinfo{journal}{\emph{Preprint arXiv:2601.22753}}
  (\bibinfo{year}{2026}).
\newblock


\bibitem[Serré et~al\mbox{.}(2025)]%
        {Serre2025}
\bibfield{author}{\bibinfo{person}{Gaëtan Serré}, \bibinfo{person}{Argyris
  Kalogeratos}, {and} \bibinfo{person}{Nicolas Vayatis}.}
  \bibinfo{year}{2025}\natexlab{}.
\newblock \showarticletitle{{Stein Boltzmann Sampling}: A Variational Approach
  for Global Optimization}. In \bibinfo{booktitle}{\emph{\International
  \Conference on Artificial Intelligence and Statistics}}.
\newblock


\bibitem[Serré et~al\mbox{.}(2026b)]%
        {Serre2026Unifying}
\bibfield{author}{\bibinfo{person}{Gaëtan Serré}, \bibinfo{person}{Argyris
  Kalogeratos}, {and} \bibinfo{person}{Nicolas Vayatis}.}
  \bibinfo{year}{2026}\natexlab{b}.
\newblock \showarticletitle{A Unifying Framework for Global Optimization: From
  Theory to Formalization}.
\newblock \bibinfo{journal}{\emph{Preprint arXiv:2508.20671}}
  (\bibinfo{year}{2026}).
\newblock


\bibitem[Tugaut(2014)]%
        {Tugaut2014}
\bibfield{author}{\bibinfo{person}{Julian Tugaut}.}
  \bibinfo{year}{2014}\natexlab{}.
\newblock \showarticletitle{Phase transitions of {McKean--Vlasov} processes in
  double-wells landscape}.
\newblock \bibinfo{journal}{\emph{Stochastics An \International \Journal of
  Probability and Stochastic Processes}} (\bibinfo{year}{2014}).
\newblock


\bibitem[Virtanen et~al\mbox{.}(2020)]%
        {Virtanen2020}
\bibfield{author}{\bibinfo{person}{Pauli Virtanen}, \bibinfo{person}{Ralf
  Gommers}, \bibinfo{person}{Travis~E. Oliphant}, \bibinfo{person}{Matt
  Haberland}, \bibinfo{person}{Tyler Reddy}, \bibinfo{person}{David
  Cournapeau}, \bibinfo{person}{Evgeni Burovski}, \bibinfo{person}{Pearu
  Peterson}, \bibinfo{person}{Warren Weckesser}, \bibinfo{person}{Jonathan
  Bright}, {and} \bibinfo{person}{et al.}} \bibinfo{year}{2020}\natexlab{}.
\newblock \showarticletitle{SciPy 1.0: fundamental algorithms for scientific
  computing in Python}.
\newblock \bibinfo{journal}{\emph{Nature Methods}} (\bibinfo{year}{2020}).
\newblock


\end{thebibliography}

\end{document}